\documentclass{article}

\usepackage[french,british]{babel}

\usepackage[T1]{fontenc}

\usepackage{amsmath}

\usepackage{amsfonts}
\usepackage{amsbsy}
\usepackage{amssymb}

\usepackage{amsthm}

\theoremstyle{plain}

\theoremstyle{definition}
\newtheorem{definition}{Definition}

\newtheorem{lemma}{Lemma}
\newtheorem{theorem}{Theorem}
\newtheorem{corollary}{Corollary}

\usepackage{latexsym}
\usepackage{amsmath}
\usepackage{amsfonts}
\usepackage{amssymb}
\usepackage{psfrag}
\usepackage{graphicx}

\begin{document}

\title{Filtering of Multidimensional Stationary Sequences with Missing Observations}

\author
{Oleksandr Masyutka\thanks
{Department of Mathematics and Theoretical Radiophysics,
Taras Shevchenko National University of Kyiv, Kyiv 01601, Ukraine, masyutkaAU@bigmir.net},
Mikhail Moklyachuk\thanks
{Department of Probability Theory, Statistics and Actuarial
Mathematics, Taras Shevchenko National University of Kyiv, Kyiv 01601, Ukraine, Moklyachuk@gmail.com},
Maria Sidei\thanks
{Department of Probability Theory, Statistics and Actuarial
Mathematics, Taras Shevchenko National University of Kyiv, Kyiv 01601, Ukraine, marysidei4@gmail.com
}
     }
\date{\today}

\maketitle

\renewcommand{\abstractname}{Abstract}
\begin{abstract}
   The  problem  of mean-square optimal linear estimation of linear functionals which depend on the unknown values of a multidimensional stationary stochastic sequence from observations of the sequence with a noise and missing observations is considered. Formulas for calculating the mean-square errors and the spectral characteristics of the optimal linear estimates of the functionals are proposed under the condition of spectral certainty, where spectral densities of the sequences are exactly known. The minimax (robust) method of estimation is applied in the case where spectral densities are not known exactly while some sets of admissible spectral densities are given. Formulas that determine the least favorable spectral densities and minimax spectral characteristics are proposed for some special sets of admissible densities.
\end{abstract}

\vspace{2ex}

\vspace{2ex}

\section{Introduction}

The problem of estimation of the unknown values of stochastic sequences and processes is of constant interest in the theory of stochastic processes. The formulation of the interpolation, extrapolation and filtering problems for stationary stochastic sequences with known spectral densities and reducing them to the corresponding problems of the theory of functions belongs to
Kolmogorov (see, for example, selected works by Kolmogorov, 1992). Effective methods of solution of the estimation problems for stationary stochastic sequences and processes were developed by Wiener (1966) and Yaglom (1987). Further results are presented in the books  by Rozanov (1967) and Hannan (1970).
The crucial assumption of most of the methods developed for estimating the unobserved values of stochastic processes is that the spectral densities of the involved stochastic processes are exactly known. However, in practice complete information on the spectral densities is impossible in most cases.
In this situation one finds parametric or nonparametric estimate of the unknown spectral density and then apply one of the traditional estimation methods provided that the selected density is the true one. This procedure can result in significant increasing of the value of error as Vastola \& Poor (1983) have demonstrated with the help of some examples.
To avoid this effect one can search the estimates which are optimal for all densities from a certain class of admissible spectral densities. These estimates are called minimax since they minimize the maximum value of the error.
The paper by Grenander (1957) was the first one where this approach to extrapolation problem for stationary processes was proposed.
Several models of spectral uncertainty and minimax-robust methods of data processing can be found in the survey paper by Kassam \& Poor (1985).
Franke (1984, 1985), Franke \& Poor ()1984 investigated the minimax extrapolation  and filtering problems for stationary sequences with the help of convex optimization methods. This approach makes it possible to find equations that determine the least favorable spectral densities for different classes of densities. In the papers by Moklyachuk (1990 -- 2015) the problems of extrapolation, interpolation and filtering  for  functionals which depend on the unknown values of stationary processes and sequences are investigated.
The estimation problems for functionals which depend on the unknown values of multivariate stationary stochastic processes is the aim of the investigation by  Moklyachuk \& Masyutka (2008 --20012). Dubovets'ka, Masyutka \& Moklyachuk (2012), Dubovets'ka \& Moklyachuk (2013 -- 2014), Moklyachuk \& Golichenko (2016) investigated the interpolation, extrapolation and filtering problems for periodically correlated stochastic sequences.
In the papers by Luz \& Moklyachuk (2012 -- 2017)  results of investigation of the estimation problems for functionals which depend on the unknown values of stochastic sequences with stationary increments are described.
Results of investigations of the prediction problem for stationary stochastic sequences with missing observations are presented in the papers by Bondon (2002, 2005), Kasahara, Pourahmadi \& Inoue (2007, 2009).
In papers by Moklyachuk \& Sidei (2015 -- 2017) results of investigations of the interpolation, extrapolation and filtering  problems for stationary stochastic sequences and processes with missing observations are proposed.

In this paper we investigate the problem of the mean-square optimal estimation of the linear functional
$A\vec{\xi}=\sum\limits_{j \in  Z ^ S}\vec{a}(j)^\top\vec{\xi}(-j)$ which depends on the unknown values of a multidimensional stationary sequence $\{\vec{\xi}(j),j\in\mathbb{Z}\}$ from the observations of the sequence $\vec{\xi}(j)+\vec{\eta}(j)$ at points $j\in\mathbb{Z_{-}}\backslash S $, where $\{\vec{\eta}(j), j\in  \mathbb{Z}\}$  is uncorrelated with $\vec{\xi}(j)$ multidimensional stationary sequence,
$S=\bigcup\limits_{l=1}^{s}\{ -(M_{l}+N_{l}),  \ldots,  -M_{l} \}$, $Z^S = \{1,2,\ldots\}\backslash S^{+}$, $S^{+}=\bigcup\limits_{l=1}^{s}\{ M_{l},  \ldots,  M_{l}+N_{l} \}$, $M_0=0$, $N_0=0$.
The problem is considered in the case where both spectral densities of the sequences $\{\vec{\xi}(j), j\in \mathbb{Z}\}$ and $\{\vec{\eta}(j), j\in  \mathbb{Z}\}$ are known. In this case we derive the spectral characteristic and the mean-square error of the optimal estimate using the method of projection in the Hilbert space of random variables with finite second moments proposed by
Kolmogorov (see, for example, selected works by Kolmogorov, 1992). In the case of spectral uncertainty, where the spectral densities of the sequences are not exactly known while a set of admissible spectral densities is given, the minimax method is applied. Formulas for determination the least favorable spectral densities and the minimax-robust spectral characteristics of the optimal estimates of the functional are proposed for some specific classes of admissible spectral densities.

\section{Hilbert space projection method of filtering}

Consider multidimensional stationary stochastic sequences  $\vec{ \xi}(j)=\left \{ \xi_ {k} (j) \right \}_{k = 1} ^ {T},\,j\in \mathbb{Z}$ and $\vec{ \eta}(j)=\left \{ \eta_ {k} (j) \right \}_{k = 1} ^ {T},\,j\in \mathbb{Z}$ with absolutely continuous spectral  functions and correlation functions of the form
\[
R_{\xi}(n)=E\vec{\xi}(j+n)(\vec{\xi}(j))^{*}=\frac{1}{2\pi}\int\limits_{-\pi}^{\pi}e^{in\lambda}F(\lambda)d\lambda,\]
\[
R_{\eta}(n)=E\vec{\eta}(j+n)(\vec{\eta}(j))^{*}=\frac{1}{2\pi}\int\limits_{-\pi}^{\pi}e^{in\lambda}G(\lambda)d\lambda,
\]
where $F(\lambda)=\left\{f_{kl}(\lambda)\right\}_{k,l=1}^T$, $G(\lambda)=\left\{g_{kl}(\lambda)\right\}_{k,l=1}^T$ are the spectral densities of the sequences  $\{\vec{\xi}(j), j\in \mathbb{Z}\}$ and $\{\vec{\eta}(j), j\in  \mathbb{Z}\}$ respectively.
\
We will suppose that the spectral densities $F(\lambda)$ and $G(\lambda)$ satisfy the minimality condition
\begin{equation}\label{minimal}
\int\limits_{-\pi}^{\pi}\left(F(\lambda)+G(\lambda)\right)^{-1}d\lambda<\infty.
\end{equation}
This condition is necessary and sufficient in order that the error-free filtering of unknown values of the sequences is impossible  (see, for example, Rozanov, 1967)

The stationary stochastic sequences  $\vec{\xi}(j)$ and $\vec{\eta}(j)$ admit the following spectral decomposition  (see, for example, Gikhman \& Skorokhod, 2004, or Karhunen, 1947)
\begin{equation} \label{ksi}
\xi(j)=\int\limits_{-\pi}^{\pi}e^{ij\lambda}Z_{\xi}(d\lambda), \hspace{1cm}
\eta(j)=\int\limits_{-\pi}^{\pi}e^{ij\lambda}Z_{\eta}(d\lambda),
\end{equation}
where $Z_{\xi}(d\lambda)$ and $Z_{\eta}(d\lambda)$ are orthogonal stochastic measures defined on $[-\pi,\pi)$ such that the following relations hold true
\[
EZ_{\xi}(\Delta_1)(Z_{\xi}(\Delta_2))^*=\frac{1}{2\pi}\int_{\Delta_1\cap\Delta_2}F(\lambda)d\lambda, \]
\[
 EZ_{\eta}(\Delta_1)(Z_{\eta}(\Delta_2))^*=\frac{1}{2\pi}\int_{\Delta_1\cap\Delta_2}G(\lambda)d\lambda.
\]

Suppose that we have observations of the sequence $\vec{\xi}(j)+\vec{\eta}(j)$ at  points $j\in\mathbb{Z_{-}}\backslash S $, where $S=\bigcup\limits_{l=1}^{s}\{ -(M_{l}+N_{l}),  \ldots,  -M_{l} \}$. The problem is to find the mean-square optimal linear estimate of the functional
$$A\vec{\xi}=\sum\limits_{j \in  Z ^ S}\vec{a}(j)^\top\vec{\xi}(-j),$$
which depends on the unknown values of the sequence $\vec{\xi}(j)$, $Z^S = \{1,2,\ldots\}\backslash S^{+}$, $S^{+}=\bigcup\limits_{l=1}^{s}\{ M_{l},  \ldots,  M_{l}+N_{l} \}$.

Suppose that the coefficients  $\{\vec{a}(j), j= 0, 1, \ldots\}$ defining the functional $A\vec{\xi}$ satisfy the following condition:
\begin{equation}\label{umovu2mom}
\sum\limits_{j \in  Z ^ S}\sum_{k=1}^T\left|a_k(j)\right|<\infty
\end{equation}
This condition ensures that the functional  $A\vec{\xi}$ has a finite second moment.

It follows from the spectral decomposition of the sequence $\vec{\xi}(j)$ that the functional $A\vec{\xi}$ can be represented in the following form
\begin{equation*}
A\vec{\xi}=\int\limits_{-\pi}^{\pi}(A(e^{i\lambda}))^\top Z_{\xi}(d\lambda), \quad
  A(e^{i\lambda})=\sum\limits_{j\in Z^S}\vec a(j)e^{-ij\lambda}.
 \end{equation*}

Consider the values $\xi_k(j)$ and $\eta_k(j)$ as the elements of the Hilbert space $H=L_2(\Omega,\mathcal{F},P)$ generated by random variables $\xi$ with zero mathematical expectations, $E\xi=0$,  finite variations, $E|\xi|^2<\infty$, and the inner product $(\xi,\eta)=E\xi\overline{\eta}$. Denote by $H^s(\xi+\eta)$ the closed linear subspace generated by elements $\{\xi_k(j)+\eta_k(j): j\in \mathbb{Z}_{-}\backslash S, k=\overline{1,T}\}$  in the Hilbert space  $H=L_2(\Omega,\mathcal{F},P)$. Denote by $L_{2} (F+G)$ the Hilbert space of vector-valued functions $ \vec{a}( \lambda )= \left \{a_{k} ( \lambda ) \right \}_{k=1}^{T} $ such that
 \[ \int_{- \pi}^{ \pi} \vec{a}( \lambda )^{ \top} \left(F(\lambda)+G(\lambda)\right) \overline{ \vec{a}( \lambda )}d \lambda < \infty. \] Denote by $L_2^s(F+G)$ the subspace of $L_2(F+G)$ generated by functions of the form
 \[e^{in \lambda} \delta_{k} , \; \delta_{k} = \left \{ \delta_{kl} \right \}_{l=1}^{T} , \; k= {1,\dots,T}, \; n \in Z_{-} \backslash S. \]

The mean square optimal linear estimate $\hat{A}\vec{\xi}$ of the functional $A\vec{\xi}$ from observations of the sequence $\vec{\xi}(j)+\vec{\eta}(j)$ can be represented in the form
 \begin {equation}\label{sp-har-general}
\hat{A}\vec{\xi}=\int\limits_{-\pi}^{\pi}(h(e^{i\lambda}))^\top(Z_{\xi}(d\lambda)+ Z_{\eta}(d\lambda)),
 \end{equation}
where $h(e^{i\lambda})=\left\{h_k(e^{i\lambda})\right\}_{k=1}^T\in L_2^s(F+G)$ is the spectral characteristic of the estimate.

The mean square error  $\Delta(h;F,G)$ of the estimate $\hat{A}\vec{\xi}$ is given by the formula
\begin{multline*}
\Delta(h;F,G)=E\left|A\vec{\xi}-\hat{A}\vec{\xi}\right|^2=
\\
=\frac{1}{2\pi}\int\limits_{-\pi}^{\pi}\left(A(e^{i\lambda})-h(e^{i\lambda})\right)^\top F(\lambda)\overline{\left(A(e^{i\lambda})-h(e^{i\lambda})\right)}d\lambda+
\\
 +\frac{1}{2\pi}\int\limits_{-\pi}^{\pi}\left(h(e^{i\lambda})\right)^\top G(\lambda)\overline{\left(h(e^{i\lambda})\right)}d\lambda.
\end{multline*}

The Hilbert space projection method proposed by A. N. Kolmogorov \cite{Kolmogorov} makes it possible to find the spectral characteristic  $h(e^{i\lambda})$ and the mean square error $\Delta(h;F,G)$ of the optimal linear estimate of the functional  $A\vec{\xi}$ in the case where   spectral densities $F(\lambda)$ and $G(\lambda)$ of the sequences are exactly known and the minimality condition (\ref{minimal}) is satisfied. According to this method the optimal estimation of the functional $A\vec{\xi}$ is a projection of the element $A\vec{\xi}$ of the space $H$  on  the space $H^s(\xi+\eta)$. It can be found from the following conditions:
\begin{equation*} \begin{split}
1)& \hat{A}\vec{\xi} \in H^s(\xi+\eta), \\
2)& A\vec{\xi}-\hat{A}\vec{\xi} \bot  H^s(\xi+\eta).
\end{split} \end{equation*}

It follows from the second condition that the spectral characteristic  $h(e^{i\lambda})$  for any $ j\in \mathbb{Z_{-}}\backslash S $ satisfies the equations
\begin{equation*}
\frac{1}{2\pi}\int\limits_{-\pi}^{\pi} \left(A(e^{i\lambda})- h(e^{i\lambda})\right)^\top F(\lambda)e^{-ij\lambda}d\lambda -\frac{1}{2\pi}\int\limits_{-\pi}^{\pi} (h(e^{-i\lambda}))^\top G(\lambda)e^{ij\lambda}d\lambda=\vec{0}.
\end{equation*}
The last relation is equivalent to equations
\begin{equation*} \begin{split}
\frac{1}{2\pi}\int\limits_{-\pi}^{\pi} \left[(A(e^{i\lambda}))^\top F(\lambda)- (h(e^{i\lambda}))^\top(F(\lambda)+G(\lambda))\right]e^{-ij\lambda}d\lambda=\vec{0},\; j\in \mathbb{Z}_{-}\backslash S.
\end{split} \end{equation*}
Hence the function $\left[(A(e^{i\lambda}))^\top F(\lambda)-(h(e^{i\lambda}))^\top(F(\lambda)+G(\lambda))\right]$
is of the form $$(A(e^{i\lambda}))^\top F(\lambda)-(h(e^{i\lambda}))^\top(F(\lambda)+G(\lambda))=(C(e^{i\lambda}))^\top,$$
where
$$C(e^{i\lambda})=\sum\limits_{j \in S}\vec{c}(j)e^{ij\lambda}+\sum\limits_{j =0}^{\infty}\vec{c}(j)e^{ij\lambda}.$$
Here $\vec{c}(j),\; j \in  U = S \cup\{0, 1, 2, \ldots\}$ are unknown coefficients that we have to find.

From the last relation we deduce that the spectral characteristic of the optimal linear estimate $\hat{A}\vec{\xi}$ is of the form
\begin{equation} \label{sphar} \begin{split}
(h(e^{i\lambda}))^\top=(A(e^{i\lambda}))^\top F(\lambda)(F(\lambda)+G(\lambda))^{-1}-(C(e^{i\lambda}))^\top(F(\lambda)+G(\lambda))^{-1}.
\end{split} \end{equation}

It follows from the first condition, $\hat{A}\vec{\xi} \in H^s(\xi+\eta)$, which determine the optimal linear estimate of the functional  $A\vec{\xi}$, that the Fourier coefficients of the function  $h(e^{i\lambda})$ are equal to zero for $k \in U$,
$$\frac{1}{2\pi}\int\limits_{-\pi}^{\pi} h(e^{i\lambda})e^{-ik\lambda }d\lambda=\vec{0}, \; k \in U,$$
namely
\begin{multline*}
\frac{1}{2\pi}\int\limits_{-\pi}^{\pi}\bigg((A(e^{i\lambda}))^\top F(\lambda)(F(\lambda)+G(\lambda))^{-1}-
\\
- (C(e^{i\lambda}))^\top(F(\lambda)+G(\lambda))^{-1}\bigg)e^{-ik\lambda}d\lambda=\vec{0}, \; k\in  U.
\end{multline*}

We will use the last equality to find equations which determine the unknown coefficients $\vec{c}(j), j\in U.$
After disclosing the brackets we get the relation
\begin{multline}\label{3}
\sum\limits_{j  \in Z^S}\vec{a}(j)^\top\frac{1}{2\pi}\int\limits_{-\pi}^{\pi}F(\lambda)(F(\lambda)+G(\lambda))^{-1}e^{-i(k+j)\lambda} d\lambda-
\\
-\sum\limits_{j \in S}\vec{c}(j)\frac{1}{2\pi}\int\limits_{-\pi}^{\pi}(F(\lambda)+G(\lambda))^{-1}e^{-i(k-j)\lambda}d\lambda-
\\
-\sum\limits_{j =0}^{\infty}\vec{c}(j)\frac{1}{2\pi}\int\limits_{-\pi}^{\pi}(F(\lambda)+G(\lambda))^{-1}e^{-i(k-j)\lambda}d\lambda=0,\; k\in U.
\end{multline}

For the functions
\[(F(\lambda)+G(\lambda))^{-1},\quad F(\lambda)(F(\lambda)+G(\lambda))^{-1},\quad F(\lambda)(F(\lambda)+G(\lambda))^{-1}G(\lambda)\]
we introduce the Fourier coefficients
\begin{equation}\label{brq1}\begin{split}
B(k,j)=\frac{1}{2\pi} \int\limits_{-\pi}^{\pi}(F(\lambda)+G(\lambda))^{-1} e^{-i(k-j)\lambda} d\lambda,\\
\end{split}\end{equation}
\begin{equation}\label{brq2}\begin{split}
R(k,j)=\frac{1}{2\pi} \int\limits_{-\pi}^{\pi}F(\lambda)(F(\lambda)+G(\lambda))^{-1} e^{-i(k+j)\lambda}d\lambda,\\
\end{split}\end{equation}
\begin{equation}\label{brq3}\begin{split}
Q(k,j)=\frac{1}{2\pi} \int\limits_{-\pi}^{\pi}F(\lambda)(F(\lambda)+G(\lambda))^{-1}G(\lambda) e^{-i(k-j)\lambda}d\lambda.
\end{split}\end{equation}

Using the introduced notations we can verify that the equality (\ref{3})  is equivalent to the following system of equations:
\begin{equation*}\begin{split}
\sum\limits_{j \in Z ^S}R(k,j)\vec{a}(j)= \sum\limits_{j \in S}B(k,j)\vec{c}(j)+
\sum\limits_{j =0}^{\infty}B(k,j)\vec{c}(j),\hspace{1cm} k\in U.
\end{split} \end{equation*}

Let us introduce notations  $\vec{a}(j)=\vec{0},\; j \in S$, $\vec{a}(0)=\vec{0}$ and $\vec{a}(j)=\vec{0},\; j \in S^{+}$. Thus, we can write
\begin{equation*}\begin{split}
\sum\limits_{j \in U}R(k,j)\vec{a}(j)=\sum\limits_{j \in S}B(k,j)\vec{c}(j)+\sum\limits_{j =0}^{\infty}B(k,j)\vec{c}(j),\hspace{1cm} k\in U.
\end{split} \end{equation*}

The last equations can be rewritten in the following form
\begin{equation}\label{rn}\begin{split}
\bold{R}\vec{\bold{a}}=\bold{B} \vec{\bold{c}},
\end{split} \end{equation}
where  $\vec{\bold{c}}$ is a vector   constructed from the unknown coefficients  $\vec{c}(j), j\in U$, vector $\vec{\bold{a}}$ has the same with the vector  $\vec{\bold{c}}$ dimension,  it is of the form  $$\vec{\bold{a}}^\top=(\vec{0}_{0}^\top, \vec{a}_1^\top, \vec{0}_1^\top, \vec{a}_2^\top,  \vec{0}_2^\top,\ldots \vec{a}_i^\top,\vec{0}_i^\top,\ldots, \vec{a}_s^\top, \vec{0}_s^\top, \vec{a}_{s+1}^\top),$$ where $\vec{0}_{0}$ is the vector which consists form    $(|S|+1)T$ zeros, where $|S|=\sum\limits_{k=1}^{s}(N_k+1)$ is the amount of missing values,
vectors $\vec{0}_i$, $i=1,2,\ldots, s,$  consist from  $(N_i + 1)T$ zeros,  vectors $$\vec{a}_1^\top=(\vec{a}(1)^\top,\ldots,\vec{a}(M_1-1)^\top),$$ $$\vec{a}_i^\top=(\vec{a}(M_{i-1}+N_{i-1}+1)^\top,\ldots,\vec{a}(M_i-1)^\top),  \quad i=2,\ldots, s,$$   $$\vec{a}_{s+1}^\top=(\vec{a}(M_s+N_s+1)^\top, \vec{a}(M_s+N_s+2)^\top,\ldots),$$  are constructed from the coefficients that determine the functional $A\vec{\xi}$.

$\bold{B}$ is a linear operator in the space $\ell_2$ which is defined by the matrix
\begin{equation*}\label{matrix2}
B= \left( \begin{array}{ccccc}
B_{s, s}& B_{s, s-1}&\ldots&B_{s, 1}&B_{s, n}\\
B_{s-1, s}& B_{s-1, s-1}&\ldots&B_{s-1, 1}&B_{s-1, n}\\
\vdots&\vdots&\ddots &\vdots&\vdots\\
B_{1, s}& B_{1, s-1}&\ldots&B_{1, 1}&B_{1, n}\\
B_{n, s}& B_{n, s-1}&\ldots&B_{n, 1}&B_{n, n}
\end{array}\right),
\end{equation*}
where   elements in the last column and the last row are compound matrices constructed from the block-matrices
\begin{equation*}\begin{split}
&B_{l,n}(k,j)= B(k,j),\; l= 1,2,\ldots,s, \; k= -M_{l}-N_l, \ldots,-M_{l}, \; j= 0, 1, 2, \ldots,\\
&B_{n,m}(k,j)= B(k,j), \; m= 1,2,\ldots,s, \; k= 0, 1, 2,\ldots, \; j= -M_{m}-N_m, \ldots,-M_{m},\\
& B_{n,n}(k,j)= B(k,j), \; k,j= 0, 1, 2,\ldots,
\end{split} \end{equation*}
and other elements of  matrix $B$ are the compound matrices with elements of the form
\begin{multline*} \label{matr32}
B_{l,m}(j,k)= b(k,j), \; l,m=1,2,\ldots,s, \\
 k= -M_{l}-N_l, \ldots,-M_{l},  j= -M_{m}-N_m, \ldots,-M_{m}.
\end{multline*}

The linear operator $\bold{R}$  in the space  $\ell_2$ is defined by the corresponding matrix in the same manner.

The unknown coefficients $\vec{c}(k),k\in U$, which are defined by the equations (\ref{rn}), can be calculated by the formula
$$\vec{c}(k)=(\bold{B}^{-1}\bold{R}\vec{\bold{a}})(k),$$
where $(\bold{B}^{-1}\bold{R}\vec{\bold{a}})(k)$ is the $k$ component of the vector $\bold{B}^{-1}\bold{R}\vec{\bold{a}}$.
The formula for calculating the spectral characteristic $h(e^{i\lambda})$ of the estimate $\hat{A}\vec{\xi}$ is of the form
\begin{equation} \label{4} \begin{split}
(h(e^{i\lambda}))^\top=(A(e^{i\lambda}))^\top F(\lambda)(F(\lambda)+G(\lambda))^{-1} - \\
-\left(\sum\limits_{k \in U}(\bold{B}^{-1}\bold{R}\vec{\bold{a}})(k)e^{ik\lambda}\right)^\top(F(\lambda)+G(\lambda))^{-1}.
\end{split} \end{equation}

The mean square error of the estimate  $\hat{A}\vec{\xi}$  can be calculated by the formula
\[
\Delta(F,G)=E\left|A\vec{\xi}-\hat{A}\vec{\xi}\right|^2=\frac{1}{2\pi}\int\limits_{-\pi}^{\pi}(r_G(\lambda))^\top F(\lambda)\overline{r_G(\lambda)}d\lambda+\]
\begin{equation} \label{55} \begin{split}
+\frac{1}{2\pi}\int\limits_{-\pi}^{\pi}(r_F(\lambda))^\top G(\lambda)\overline{r_F(\lambda)}d\lambda
&=\langle\bold{R}\vec{\bold{a}},\bold{B}^{-1}\bold{R}\vec{\bold{a}}\rangle+\langle\bold{Q}\vec{\bold{a}},\vec{\bold{a}}\rangle,
\end{split} \end{equation}
where
\[(r_F(\lambda))^\top=\left((A(e^{i\lambda}))^\top F(\lambda)-\left(\sum\limits_{k \in U}(\bold{B}^{-1}\bold{R}\vec{\bold{a}})(k) e^{ik\lambda}\right)^\top\right)(F(\lambda)+G(\lambda))^{-1},\]
\[(r_G(\lambda))^\top=\left((A(e^{i\lambda}))^\top G(\lambda)+\left(\sum\limits_{k \in U}(\bold{B}^{-1}\bold{R}\vec{\bold{a}})(k) e^{ik\lambda}\right)^\top\right)(F(\lambda)+G(\lambda))^{-1},\]
and $\bold{Q}$ is the linear operator in the space $\ell_2$ defined by matrix with coefficients $Q(k,j)$,  $k, j \in U$ in the same way as the operator $\bold{B}$ is defined.

Let us summarize results and present them in the form of a theorem.
\begin{theorem}{Theorem 1.}{}
Let $\vec{\xi}(j)$ and  $\vec{\eta}(j)$ be uncorrelated multidimensional stationary sequences with spectral densities $F(\lambda)$ and $G(\lambda)$ which satisfy the minimality condition (\ref{minimal}). The spectral characteristic   $h(e^{i\lambda})$ and the mean square error  $\Delta(F,G)$ of the optimal linear estimate of the functional  $A\vec{\xi}$ which depends on the unknown values of the sequence  $\vec{\xi}(j)$ based on observations of the sequence  $\vec{\xi}(j)+\vec{\eta}(j),$ $j\in \mathbb{Z}_{-}\backslash S$ can be calculated by formulas (\ref{4}), (\ref{55}).
\end{theorem}

Consider the problem of the mean-square optimal linear estimation of the functional
$$A\vec{\xi}=\sum\limits_{j \in Z^S}\vec{a}(j)^\top\vec{\xi}(-j),$$
which depends on the unknown values of the sequence  $\vec{\xi}(j)$ from  observations of the sequence  $\vec{\xi}(j)+\vec{\eta}(j)$ at points $j\in\mathbb{Z_{-}}\backslash S $,  $S=\{ -(M+N),  \ldots,  -M \}$, $Z^S = \{1,2,\ldots\}\backslash S^{+}$, $S^{+}=\{ M,  \ldots,  M+N \}$.

From Theorem 1 the following corollary can be derived for this problem.
\begin{corollary}{Corollary 1.}{}%
Let $\vec{\xi}(j)$ and  $\vec{\eta}(j)$ be uncorrelated multidimensional stationary sequences with spectral densities $F(\lambda)$ and $G(\lambda)$ which satisfy the minimality condition (\ref{minimal}). The spectral characteristic   $h(e^{i\lambda})$ and the mean square error  $\Delta(F,G)$ of the optimal linear estimate of the functional  $A\vec{\xi}$ which depends on the unknown values of the sequence  $\vec{\xi}(j)$ based on observations of the sequence  $\vec{\xi}(j)+\vec{\eta}(j),$ $j\in \mathbb{Z}_{-}\backslash S$ can be calculated by formulas (\ref{4ss}), (\ref{55ss})
\begin{multline}\label{4ss}
(h(e^{i\lambda}))^\top=(A(e^{i\lambda}))^\top F(\lambda)(F(\lambda)+G(\lambda))^{-1}-
\\
-\left(\sum\limits_{k \in U}(\bold{B}^{-1}\bold{R}\vec{\bold{a}})_ke^{ik\lambda}\right)^\top(F(\lambda)+G(\lambda))^{-1},
\end{multline}
\begin{equation} \label{55ss} \begin{split}
\Delta(h;F,G)=\langle\bold{R}\vec{\bold{a}},\bold{B}^{-1}\bold{R}\vec{\bold{a}}\rangle+
\langle\bold{Q}\vec{\bold{a}},\vec{\bold{a}}\rangle,
\end{split} \end{equation}
$\bold{B}$, $\bold{R}$, $\bold{Q}$ are linear operators in the space $\ell_2$ defined by compound matrices with coefficients $B(k,j)$, $R(k,j)$, $Q(k,j)$,  $k, j \in U$, $(U=S\cup\{0,1,2,\ldots\})$.
For example, the matrix $B$ is of the form
\begin{equation*} \label{matrix2}
B= \left( \begin{array}{cc}
B_{s, s}&B_{s, n}\\
B_{n, s}&B_{n, n}
\end{array}\right),
\end{equation*}
where its components are matrices constructed from the block-matrices
 \begin{equation*}\begin{split}
&B_{s,n}(k,j)= B(k,j), \hspace{0.5cm} k= -M -N , \ldots,-M , \quad j= 0, 1, 2, \ldots,\\
&B_{n,s}(k,j)= B(k,j), \hspace{0.5cm} k= 0, 1, 2,\ldots,\quad j= -M -N , \ldots,-M ,\\
&B_{n,n}(k,j)= B(k,j), \hspace{0.5cm} k,j= 0, 1, 2,\ldots,\\
&B_{s,s}(k,j)= B(k,j), \hspace{0.5cm} k= -M -N , \ldots,-M , \quad j= -M -N , \ldots,-M .
\end{split}\end{equation*}
\end{corollary}

Consider the problem of the  mean-square optimal linear estimation of the functional
$$A\vec{\xi}=\sum\limits_{j \in Z^S}\vec{a}(j)^\top\vec{\xi}(-j),$$
which depends on the unknown values of the sequence  $\vec{\xi}(j)$ from  observations of the sequence  $\vec{\xi}(j)+\vec{\eta}(j)$ at points $j\in\mathbb{Z_{-}}\backslash \{-s\} $,  $Z^S = \{1,2,\ldots\} \backslash \{s\}$.

It follows from Theorem 1 that the following corollary  holds true.
\begin{corollary}{Corollary 2.}{}%
Let $\vec{\xi}(j)$ and $\vec{\eta}(j)$ be uncorrelated multidimensional stationary sequences with spectral densities $F(\lambda)$ and $G(\lambda)$ which satisfy the minimality condition (\ref{minimal}). The spectral characteristic   $h(e^{i\lambda})$ and the mean square error  $\Delta(F,G)$ of the optimal linear estimate of the functional  $A\vec{\xi}$ which depends on the unknown values of the sequence  $\vec{\xi}(j)$ based on observations of the sequence  $\vec{\xi}(j)+\vec{\eta}(j),$ $j\in \mathbb{Z}_{-}\backslash \{-s\}$ can be calculated by formulas  (\ref{4s}), (\ref{55s})
\begin{multline}\label{4s}
(h(e^{i\lambda}))^\top=(A(e^{i\lambda}))^\top F(\lambda)(F(\lambda)+G(\lambda))^{-1}-\\
-\left(\sum\limits_{k \in U}(\bold{B}^{-1}\bold{R}\vec{\bold{a}})_ke^{ik\lambda}\right)^\top(F(\lambda)+G(\lambda))^{-1},
\end{multline}
\begin{equation} \label{55s} \begin{split}
\Delta(h;F,G)=\langle\bold{R}\vec{\bold{a}},\bold{B}^{-1}\bold{R}\vec{\bold{a}}\rangle+\langle\bold{Q}\vec{\bold{a}},\vec{\bold{a}}\rangle,
\end{split} \end{equation}
$\bold{B}$, $\bold{R}$, $\bold{Q}$ are linear operators in the space $\ell_2$ defined by compound matrices with coefficients $B(k,j)$, $R(k,j)$, $Q(k,j)$,  $k, j \in U$, $(U=S\cup\{0,1,2,\ldots\})$,
\begin{equation*} \label{matrix2}
B= \left( \begin{array}{cc}
B(-s, -s)&B_{-s, n}\\
B_{n, -s}&B_{n, n}
\end{array}\right),
\end{equation*}
where elements in the last column and the last row are the matrices with the elements \begin{equation*}\begin{split}
&B_{-s,n}(k,j)= B(k,j), \hspace{0.5cm} k= -s , \quad j= 0, 1, 2, \ldots,\\
&B_{n,-s}(k,j)=B(k,j), \hspace{0.5cm} k= 0, 1, 2,\ldots, \quad j= -s,\\
&B_{n,n}(k,j)= B(k,j), \hspace{0.5cm} k,j= 0, 1, 2,\ldots.
\end{split}\end{equation*}
\end{corollary}
Consider the problem of the  mean-square optimal linear estimation of the functional
$$A_N\vec{\xi}=\sum\limits_{j \in Z^S \cap \{0, \ldots, N \}}\vec{a}(j)^\top\vec{\xi}(-j),$$
which depends on the unknown values of the sequence  $\vec{\xi}(j)$ from  observations of the sequence  $\vec{\xi}(j)+\vec{\eta}(j)$ at points $j\in \mathbb{Z}_{-}\backslash S$ where $S$ is defined in the introduction. The linear estimate of the functional  $A_N\vec{\xi}$ has the representation
\begin{equation*}
\hat{A}_N\vec{\xi}=\int\limits_{-\pi}^{\pi}(h_N(e^{i\lambda})^\top(Z_{\xi}(d\lambda)+ Z_{\eta}(d\lambda)).
 \end{equation*}

Define the vector $\vec{\bold{a}}_N$ as follows: elements with indices from the set $U\cap (S\cup\{0, \ldots, N \})$ coincide with the elements of the vector $\vec{\bold{a}}$ with the same indices and elements with indices from the set  $U\backslash(S\cup \{0, \ldots, N \})$ are zeros. $\bold{B}$, $\bold{R}$, $\bold{Q}$ are linear operators in the space $\ell_2$  defined in the Theorem  1.

The spectral characteristic   $h_N(e^{i\lambda})$ and the mean square error  $\Delta(h_N;F,G)$ of the optimal linear estimate of the functional  $A_N\vec{\xi}$ can be calculated by  formulas (\ref{4d}), (\ref{55d})
\begin{multline}\label{4d}
(h_N(e^{i\lambda}))^\top=(A_N(e^{i\lambda}))^\top F(\lambda)(F(\lambda)+G(\lambda))^{-1}-\\
-\left(\sum\limits_{k \in U}(\bold{B}^{-1}\bold{R}\vec{\bold{a}}_N)(k)e^{ik\lambda}\right)^\top(F(\lambda)+G(\lambda))^{-1},
\end{multline}
\begin{equation} \label{55d} \begin{split}
\Delta(h_N;F,G)=\langle\bold{R}\vec{\bold{a}}_N,\bold{B}^{-1}\bold{R}\vec{\bold{a}}_N\rangle+\langle\bold{Q}\vec{\bold{a}}_N,\vec{\bold{a}}_N\rangle,
\end{split} \end{equation}
 where $ A_N(e^{i\lambda})=\sum\limits_{j \in Z^S \cap \{0, \ldots, N \}} \vec{a}(j)e^{-ij\lambda}.$

The following corollary holds true.

\begin{corollary}{Corollary 3.}{}%
Let $\vec{\xi}(j)$ and  $\vec{\eta}(j)$ be multidimensional uncorrelated stationary sequences with spectral densities $F(\lambda)$ and $G(\lambda)$ which satisfy the minimality condition (\ref{minimal}). The spectral characteristic   $h_N(e^{i\lambda})$ and the mean square error  $\Delta(h_N;F,G)$ of the optimal linear estimate of the functional  $A_N\vec{\xi}$ which depends on the unknown values of the sequence $\vec{\xi}(j)$ from observation of the sequence $\vec{\xi}(j)+\vec{\eta}(j)$ at points of time $j\in \mathbb{Z}_{-}\backslash S $ can be calculated by formulas (\ref{4d}), (\ref{55d}).
\end{corollary}

\section{Minimax-robust method of filtering}

Theorem 1 and its corollaries can be applied to  filtering of the functional in the cases when spectral densities of the sequences are exactly known. If complete information on the spectral densities is impossible but the class of admissible densities is given, it is reasonable to apply the minimax-robust method of filtering which consists in minimizing the value of the mean-square error for all spectral densities from the given class. For description of minimax method we propose the following definitions ( see Moklyachuk, 2000).

 \begin{definition}{Definition 1.}{}%
 For a given class of spectral densities $D=D_F \times D_G$ the spectral densities  $F^0(\lambda) \in D_F$, $G^0(\lambda) \in D_G$ are called least favorable in the class $D$ for the optimal linear filtering of the functional $A\vec{\xi}$ if the following relation holds true
$$\Delta\left(F^0,G^0\right)=\Delta\left(h\left(F^0,G^0\right);F^0,G^0\right)=\max\limits_{(F,G)\in D_F\times D_G}\Delta\left(h\left(F,G\right);F,G\right).$$
\end{definition}

 \begin{definition}{Definition 2.}{}%
For a given class of spectral densities $D=D_F \times D_G$ the spectral characteristic $h^0(e^{i\lambda})$ of the optimal linear estimate  of the functional  $A\vec{\xi}$ is called minimax-robust if there are satisfied conditions
$$h^0(e^{i\lambda})\in H_D= \bigcap\limits_{(F,G)\in D_F\times D_G} L_2^s(F+G),$$
$$\min\limits_{h\in H_D}\max\limits_{(F,G)\in D}\Delta\left(h;F,G\right)=\max\limits_{(F,G)\in D}\Delta\left(h^0;F,G\right).$$
\end{definition}

From the introduced definitions and formulas derived above we can obtain the following statement.

\begin{lemma}{Lemma 1.}{}%
Spectral densities $F^0(\lambda)\in D_F,$ $G^0(\lambda) \in D_G$ satisfying the minimality condition (\ref{minimal}) are the least favorable in the class  $D=D_F\times D_G$ for the optimal linear filtering of the functional $A\vec{\xi}$ if operators $B^0, R^0, Q^0$ determined by the Fourier coefficients of the functions $$(F^0(\lambda)+G^0(\lambda))^{-1}, \, F^0(\lambda)(F^0(\lambda)+G^0(\lambda))^{-1}, \, F^0(\lambda)(F^0(\lambda)+G^0(\lambda))^{-1}G^0(\lambda)$$ determine a solution to the constrain optimization problem
\begin{equation} \label{extrem} \begin{split}
\max\limits_{(F,G)\in D_F\times D_G}\langle\bold{R}\vec{\bold{a}},\bold{B}^{-1}\bold{R}\vec{\bold{a}}\rangle +\langle\bold{Q}\vec{\bold{a}},\vec{\bold{a}}\rangle=\langle\bold{R}^0\vec{\bold{a}},(\bold{B}^0)^{-1}\bold{R}^0\vec{\bold{a}}\rangle+\langle\bold{Q}^0\vec{\bold{a}},\vec{\bold{a}}\rangle.
\end{split}\end{equation}
The minimax spectral characteristic $h^0=h(F^0,G^0)$ is determined by the formula (\ref{4}) if $h(F^0,G^0) \in H_D.$
\end{lemma}

The least favorable spectral densities $F^0(\lambda)$, $G^0(\lambda)$ and the minimax spectral characteristic $h^0=h(F^0,G^0)$ form a saddle point of the function $\Delta \left(h;F,G\right)$ on the set  $H_D\times D.$ The saddle point inequalities
$$\Delta\left(h;F^0,G^0\right)\geq\Delta\left(h^0;F^0,G^0\right)\geq \Delta\left(h^0;F,G\right) $$  $$ \forall h \in H_D, \forall F \in D_F, \forall G \in D_G$$
hold true if  $h^0=h(F^0,G^0)$ and  $h(F^0,G^0)\in H_D,$  where $(F^0,G^0)$ is a solution to the constrained optimization problem
\begin{equation} \label{7}
\sup\limits_{(F,G)\in D_F\times D_G}\Delta\left(h(F^0,G^0);F,G\right)=\Delta\left(h(F^0,G^0);F^0,G^0\right),
\end{equation}
where the functional $\Delta\left(h(F^0,G^0);F,G\right)$ is calculated by the formula
\begin{multline*}
\Delta\left(h(F^0,G^0);F,G\right)=\frac{1}{2\pi}\int\limits_{-\pi}^{\pi}(r_G^0(\lambda))^\top F(\lambda)\overline{r_G^0(\lambda)}d\lambda+
\\
+\frac{1}{2\pi}\int\limits_{-\pi}^{\pi}(r_F^0(\lambda))^\top G(\lambda)\overline{r_F^0(\lambda)}d\lambda,
\end{multline*}

\begin{multline*}
(r_F^0(\lambda))^\top=
\\
=\left((A(e^{i\lambda}))^\top F^0(\lambda)-\left(\sum\limits_{k \in U}((\bold{B}^0)^{-1}\bold{R}^0\vec{\bold{a}})(k) e^{ik\lambda}\right)^\top\right)(F^0(\lambda)+G^0(\lambda))^{-1},
\end{multline*}

\begin{multline*}
(r_G^0(\lambda))^\top=
\\
=\left((A(e^{i\lambda}))^\top G^0(\lambda)+\left(\sum\limits_{k \in U}((\bold{B}^0)^{-1}\bold{R}^0\vec{\bold{a}})(k) e^{ik\lambda}\right)^\top\right)(F^0(\lambda)+G^0(\lambda))^{-1}.
\end{multline*}
The constrained optimization problem  (\ref{7}) is equivalent to the unconstrained optimization problem (see, for example, Pshenichnyj, 1971)
\begin{equation} \label{8}
\Delta_D(F,G)=-\Delta(h(F^0,G^0);F,G)+\delta((F,G)\left|D_F\times D_G\right.)\rightarrow \inf,
\end{equation}
where $\delta((F,G)\left|D_F\times D_G\right.)$ is the indicator function of the set $D=D_F\times D_G$. Solution of the problem (\ref{8}) is characterized by the condition $0 \in \partial\Delta_D(F^0,G^0),$ where $ \partial\Delta_D(F^0,G^0)$ is the subdifferential of the convex functional $\Delta_D(F,G)$ at point $(F^0,G^0)$. This condition makes it possible to find the least favourable spectral densities in some special classes of spectral densities $D$ (see books by Ioffe \& Tihomirov, 1979,  Pshenichnyj, 1971, Rockafellar, 1997).

Note, that the form of the functional $\Delta\left(h^0;F,G\right)$ is convenient for application the Lagrange method of indefinite multipliers for finding solution to the problem (\ref{8}). Making use the method of Lagrange multipliers and the form of
subdifferentials of the indicator functions we describe relations that determine least favourable spectral densities in some special classes of spectral densities (see books by Moklyachuk,2008, Moklyachuk \& Masyutka, 2012 for additional details).

\section{Least favorable spectral densities in the class $D=D_0 \times D_{2\delta} $}

Consider the problem of filtering of the functional $A\vec{\xi}$ in the case where spectral densities $F(\lambda)$, $G(\lambda)$ belong to the set of admissible spectral densities $D_0 \times D_{2\delta} $
$$ D_{0}^{1} =\bigg\{F(\lambda )\left|\frac{1}{2\pi }
\int _{-\pi }^{\pi}{\rm{Tr}}\, F(\lambda )d\lambda =p\right.\bigg\},$$
$$D_{2\delta}^{1}=\left\{G(\lambda )\biggl|\frac{1}{2\pi } \int_{-\pi}^{\pi}\left|{\rm{Tr}}(G(\lambda )-G_{1} (\lambda))\right|^{2} d\lambda \le \delta\right\};$$
$$D_{0}^{2} =\bigg\{F(\lambda )\left|\frac{1}{2\pi }
\int _{-\pi}^{\pi}f_{kk} (\lambda )d\lambda =p_{k}, k=\overline{1,T}\right.\bigg\},$$
$$D_{2\delta}^{2}=\left\{G(\lambda )\biggl|\frac{1}{2\pi } \int_{-\pi}^{\pi}\left|g_{kk} (\lambda )-g_{kk}^{1} (\lambda)\right|^{2} d\lambda  \le \delta_{k}, k=\overline{1,T}\right\};$$
$$D_{0}^{3} =\bigg\{F(\lambda )\left|\frac{1}{2\pi} \int _{-\pi}^{\pi}\left\langle B_{1} ,F(\lambda )\right\rangle d\lambda  =p\right.\bigg\},$$
$$D_{2\delta}^{3}=\left\{G(\lambda )\biggl|\frac{1}{2\pi } \int_{-\pi}^{\pi}\left|\left\langle B_{2} ,G(\lambda )-G_{1}(\lambda )\right\rangle \right|^{2} d\lambda  \le \delta\right\};$$
$$D_{0}^{4} =\bigg\{F(\lambda )\left|\frac{1}{2\pi} \int
_{-\pi}^{\pi}F(\lambda )d\lambda  =P\right.\bigg\},$$
$$D_{2\delta}^{4}=\left\{G(\lambda )\biggl|\frac{1}{2\pi } \int_{-\pi}^{\pi}\left|g_{ij} (\lambda )-g_{ij}^{1} (\lambda)\right|^{2} d\lambda  \le \delta_{ij}, i,j=\overline{1,T}\right\},$$
where spectral density $G_{1} (\lambda))$ is known and fixed, $p, \delta, p_k, \delta_k, k=\overline{1,T}$, $\delta_{ij}, i,j=\overline{1,T}$, are given numbers, $P, B_1, B_2$ are given positive-definite Hermitian matrices.
From the condition $0\in \partial \Delta _{D} (F^{0} ,G^{0} )$ we find the following equations which determine the least favourable spectral densities for these given sets of admissible spectral densities.

For the first pair $D_{0}^{1}\times D_{2\delta}^{1}$ we have equations
\begin{equation} \label{eq_4_1}
(r_G^0(\lambda))^{*}(r_G^0(\lambda))^\top=\alpha^{2} (F^{0} (\lambda )+G^{0} (\lambda
))^{2} ,
\end{equation}
\begin{equation} \label{eq_4_2}
(r_F^0(\lambda))^{*}(r_F^0(\lambda))^\top=\beta ^{2} {\mathrm{Tr}}\, (G^{0}(\lambda )-G_{1} (\lambda ))(F^{0} (\lambda )+G^{0} (\lambda ))^{2},
\end{equation}
\begin{equation} \label{eq_4_3}
\frac{1}{2\pi } \int _{-\pi}^{\pi}\left|{\mathrm{Tr}}\, (G(\lambda )-G_{1} (\lambda ))\right|^{2} d\lambda  =\delta.
\end{equation}

For the second pair $D_{0}^{2}\times D_{2\delta}^{2}$ we have equations
\begin{equation}  \label{eq_4_4}
(r_G^0(\lambda))^{*}(r_G^0(\lambda))^\top=(F^{0} (\lambda )+G^{0} (\lambda
))\left\{\alpha _{k}^{2} \delta _{kl} \right\}_{k,l=1}^{T} (F^{0}
(\lambda )+G^{0} (\lambda )),
\end{equation}
\begin{multline} \label{eq_4_5}
(r_F^0(\lambda))^{*}(r_F^0(\lambda))^\top=\\
=(F^{0} (\lambda )+G^{0}(\lambda ))\left\{\beta _{k}^{2} (g_{kk}^{0} (\lambda )-g_{kk}^{1}
(\lambda ))\delta _{kl} \right\}_{k,l=1}^{T} (F^{0} (\lambda )+G^{0}(\lambda )),
\end{multline}
\begin{equation} \label{eq_4_6}
\frac{1}{2\pi } \int _{-\pi}^{\pi}\left|g_{kk} (\lambda )-g_{kk}^{1} (\lambda )\right|^{2} d\lambda  =\delta _{k},\; k=\overline{1,T}.
\end{equation}

For the third pair $D_{0}^{3}\times D_{2\delta}^{3}$ we have equations
\begin{equation}  \label{eq_4_7}
(r_G^0(\lambda))^{*}(r_G^0(\lambda))^\top=\alpha^{2} (F^{0} (\lambda )+G^{0} (\lambda
))B_{1}^{\top} (F^{0} (\lambda )+G^{0} (\lambda )),
\end{equation}
\begin{equation} \label{eq_4_8}
(r_F^0(\lambda))^{*}(r_F^0(\lambda))^\top=\beta ^{2} \left\langle B_{2},G^{0} (\lambda )-G_{1} (\lambda )\right\rangle(F^{0} (\lambda )+G^{0} (\lambda ))^2,
\end{equation}
\begin{equation} \label{eq_4_9}
\frac{1}{2\pi } \int _{-\pi}^{\pi}\left|\left\langle B_{2} ,G(\lambda )-G_{1} (\lambda )\right\rangle \right|^2d\lambda  =\delta.
\end{equation}

For the fourth pair $D_{0}^{4}\times D_{2\delta}^{4}$ we have equations
\begin{equation} \label{eq_4_10}
(r_G^0(\lambda))^{*}(r_G^0(\lambda))^\top=(F^{0} (\lambda )+G^{0} (\lambda
))\vec{\alpha}\cdot \vec{\alpha}^{*}(F^{0} (\lambda )+G^{0}
(\lambda )),
\end{equation}
\begin{multline} \label{eq_4_11}
(r_F^0(\lambda))^{*}(r_F^0(\lambda))^\top=\\=(F^{0} (\lambda )+G^{0} (\lambda ))\left\{\beta _{ij}^{} (g_{ij}^{0} (\lambda )-g_{ij}^{1}
(\lambda ))\right\}_{i,j=1}^{T} (F^{0} (\lambda )+G^{0} (\lambda )),
\end{multline}
\begin{equation} \label{eq_4_12}
\frac{1}{2\pi } \int _{-\pi}^{\pi}\left|g_{ij} (\lambda )-g_{ij}^{1} (\lambda )\right|^{2} d\lambda  =\delta_{ij},\; i,j=\overline{1,T}.
\end{equation}

The following theorem and corollaries hold true.
\begin{theorem}{Theorem 3.}{}%
Let the minimality condition (\ref{minimal}) hold true. The least favorable spectral densities  $F^0(\lambda)$, $G^0(\lambda)$  in the classes $D_0 \times D_{2\delta}$ for the optimal linear filtering of the functional $A\vec{\xi}$ are determined by relations
(\ref{eq_4_1}) -- (\ref{eq_4_3}) for the first pair $D_{0}^{1}\times D_{2\delta}^{1}$ of sets of admissible spectral densities;
(\ref{eq_4_4}) -- (\ref{eq_4_6}) for the second pair $D_{0}^{2}\times D_{2\delta}^{2}$ of sets of admissible spectral densities;
(\ref{eq_4_7}) -- (\ref{eq_4_9}) for the third pair $D_{0}^{3}\times D_{2\delta}^{3}$ of sets of admissible spectral densities;
(\ref{eq_4_10}) -- (\ref{eq_4_12}) for the fourth pair $D_{0}^{4}\times D_{2\delta}^{4}$of sets of admissible spectral densities;
constrained optimization problem (\ref{extrem}) and restrictions  on densities from the corresponding classes $D_0 \times D_{2\delta}$.  The minimax-robust spectral characteristic of the optimal estimate of the functional $A\vec{\xi}$ is determined by the formula (\ref{4}).
\end{theorem}

\begin{corollary}{Corollary 4.}{}%
Assume that the spectral density $G(\lambda)$ is known. Let the function $F^0(\lambda)+G(\lambda)$ satisfy the minimality condition (\ref{minimal}). The spectral density $F^0(\lambda)$ is the least favorable in the classes $D_0^k$, $k=\overline{1,4}$, for the optimal linear filtering of the functional $A\vec{\xi}$ if it satisfies relations (\ref{eq_4_1}), (\ref{eq_4_4}), (\ref{eq_4_7}), (\ref{eq_4_10}), respectively, and the pair $(F^0(\lambda), G(\lambda))$ is a solution of the optimization problem  (\ref{extrem}). The minimax-robust spectral characteristic of the optimal estimate of the functional $A\vec{\xi}$ is determined by formula (\ref{4}).
\end{corollary}

\begin{corollary}{Corollary 5.}{}%
Assume that the spectral density $F(\lambda)$ is known. Let the function $F(\lambda)+G^0(\lambda)$ satisfy the minimality condition (\ref{minimal}). The spectral density $G^0(\lambda)$ is the least favorable in the classes $D_{2\delta}^k$, $k=\overline{1,4}$, for the optimal linear filtering of the functional $A\vec{\xi}$ if it satisfies relations (\ref{eq_4_2}) -- (\ref{eq_4_3}), (\ref{eq_4_5}) -- (\ref{eq_4_6}), (\ref{eq_4_8}) -- (\ref{eq_4_9}), (\ref{eq_4_11}) -- (\ref{eq_4_12}), respectively, and the pair $(F(\lambda), G^0(\lambda))$ is a solution of the optimization problem  (\ref{extrem}). The minimax-robust spectral characteristic of the optimal estimate of the functional $A\vec{\xi}$ is determined by formula (\ref{4}).
\end{corollary}

\section{Least favorable spectral densities in the class $D=D_{1\delta} \times D_V^U$}

Consider the problem of filtering of the functional $A\vec{\xi}$ in the case where spectral densities $F(\lambda)$, $G(\lambda)$ belong to the set of admissible spectral densities $D_{1\delta} \times D_V^U$

$$D_{1\delta}^{1}=\left\{F(\lambda )\biggl|\frac{1}{2\pi} \int_{-\pi}^{\pi}\left|{\rm{Tr}}(F(\lambda )-F_{1} (\lambda))\right|d\lambda \le \delta\right\},$$
$${D_{V}^{U}} ^{1}  =\bigg\{G(\lambda )\bigg|{\mathrm{Tr}}\, V(\lambda
)\le {\mathrm{Tr}}\, G(\lambda )\le {\mathrm{Tr}}\, U(\lambda ), \frac{1}{2\pi } \int _{-\pi}^{\pi}{\mathrm{Tr}}\,  G(\lambda)d\lambda  =q \bigg\},$$
$$D_{1\delta}^{2}=\left\{F(\lambda )\biggl|\frac{1}{2\pi } \int_{-\pi}^{\pi}\left|f_{kk} (\lambda )-f_{kk}^{1} (\lambda)\right|d\lambda  \le \delta_{k}, k=\overline{1,T}\right\},$$
$${D_{V}^{U}} ^{2}  =\bigg\{G(\lambda )\bigg|v_{kk} (\lambda )  \le
g_{kk} (\lambda )\le u_{kk} (\lambda ), \frac{1}{2\pi} \int _{-\pi}^{\pi}g_{kk} (\lambda
)d\lambda  =q_{k} , k=\overline{1,T}\bigg\},$$
$$D_{1\delta}^{3}=\left\{F(\lambda )\biggl|\frac{1}{2\pi } \int_{-\pi}^{\pi}\left|\left\langle B_{1} ,F(\lambda )-F_{1}(\lambda )\right\rangle \right|d\lambda  \le \delta\right\},$$
$${D_{V}^{U}} ^{3}  =\bigg\{G(\lambda )\bigg|\left\langle B_{2}
,V(\lambda )\right\rangle \le \left\langle B_{2} ,G(\lambda
)\right\rangle \le \left\langle B_{2} ,U(\lambda)\right\rangle,\frac{1}{2\pi }
\int _{-\pi}^{\pi}\left\langle B_{2},G(\lambda)\right\rangle d\lambda  =q\bigg\},$$
$$D_{1\delta}^{4}=\left\{F(\lambda )\biggl|\frac{1}{2\pi} \int_{-\pi}^{\pi}\left|f_{ij} (\lambda )-f_{ij}^{1} (\lambda)\right|d\lambda  \le \delta_{ij}, i,j=\overline{1,T}\right\},$$
$${D_{V}^{U}} ^{4}=\left\{G(\lambda )\bigg|V(\lambda )\le G(\lambda
)\le U(\lambda ), \frac{1}{2\pi } \int _{-\pi}^{\pi}G(\lambda )d\lambda=Q\right\},$$
where spectral densities $F_{1} ( \lambda ), V( \lambda ),U( \lambda )$ are known and fixed, $\delta, q, \delta_k, q_k, k=\overline{1,T}$, $\delta_{ij}, i,j=\overline{1,T}$, are given numbers, $Q, B_1, B_2$ are given positive-definite Hermitian matrices.

From the condition $0\in \partial \Delta _{D} (F^{0} ,G^{0} )$ we find the following equations which determine the least favourable spectral densities for these given sets of admissible spectral densities.

For the first pair $D_{1\delta}^{1}\times {D_{V}^{U}} ^{1}$ we have equations
\begin{equation} \label{eq_5_1}
(r_G^0(\lambda))^{*}(r_G^0(\lambda))^\top=\alpha^{2}\gamma(\lambda )(F^{0}(\lambda )+G^{0}(\lambda ))^{2},
\end{equation}
\begin{equation} \label{eq_5_2}
\frac{1}{2 \pi} \int_{-\pi}^{\pi} \left|{\mathrm{Tr}}\, (F^0( \lambda )-F_{1}(\lambda )) \right|d\lambda =\delta,
\end{equation}
\begin{equation} \label{eq_5_3}
(r_F^0(\lambda))^{*}(r_F^0(\lambda))^\top=(\beta^{2} +\gamma _{1}
(\lambda )+\gamma _{2} (\lambda ))(F^{0} (\lambda )+G^{0} (\lambda))^{2},
\end{equation}
where $\left| \gamma( \lambda ) \right| \le 1$ and
\[\gamma(\lambda )={\mathrm{sign}}\; ({\mathrm{Tr}}\, (F^{0} ( \lambda )-F_{1} ( \lambda ))): \; {\mathrm{Tr}}\, (F^{0} ( \lambda )-F_{1} ( \lambda )) \ne 0,\]
$\gamma _{1} (\lambda )\le 0$ and $\gamma _{1} (\lambda )=0$ if ${\mathrm{Tr}}\,
G^{0} (\lambda )> {\mathrm{Tr}}\,  V(\lambda ),$ $\gamma _{2} (\lambda )\ge 0$ and $\gamma _{2} (\lambda )=0$ if $ {\mathrm{Tr}}\,G^{0}(\lambda )< {\mathrm{Tr}}\,  U(\lambda).$

For the second pair $D_{1\delta}^{2}\times {D_{V}^{U}} ^{2}$ we have equations
\begin{equation}\label{eq_5_4}
(r_G^0(\lambda))^{*}(r_G^0(\lambda))^\top=(F^{0}(\lambda )+G^{0} (\lambda ))\left\{\alpha_{k}^{2} \gamma_{k}(\lambda )\delta _{kl} \right\}_{k,l=1}^{T} (F^{0} (\lambda)+G^{0} (\lambda )),
\end{equation}
\begin{equation} \label{eq_5_5}
\frac{1}{2 \pi} \int_{-\pi}^{\pi} \left|f^0_{kk} ( \lambda)-f_{kk}^{1} ( \lambda ) \right| d\lambda =\delta_{k}, \; k=\overline{1,T},
\end{equation}
\begin{multline}\label{eq_5_6}
(r_F^0(\lambda))^{*}(r_F^0(\lambda))^\top=\\
=(F^{0} (\lambda)+G^{0} (\lambda ))\left\{(\beta_{k}^{2} +\gamma _{1k} (\lambda )+\gamma _{2k}(\lambda ))\delta _{kl} \right\}_{k,l=1}^{T} (F^{0} (\lambda )+G^{0}(\lambda )),
\end{multline}
where $\left| \gamma_{k} ( \lambda ) \right| \le 1$ and
\[\gamma_{k}( \lambda )={ \mathrm{sign}}\;(f_{kk}^{0}( \lambda)-f_{kk}^{1} ( \lambda )): \; f_{kk}^{0} ( \lambda )-f_{kk}^{1}(\lambda ) \ne 0, \; k= \overline{1,T},\]
$\gamma _{1k} (\lambda )\le 0$ and $\gamma _{1k} (\lambda )=0$ if $g_{kk}^{0} (\lambda )>v_{kk} (\lambda ),$ $\gamma _{2k} (\lambda )\ge 0$ and $\gamma _{2k} (\lambda )=0$ if $g_{kk}^{0} (\lambda )<u_{kk} (\lambda).$

For the third pair $D_{1\delta}^{3}\times {D_{V}^{U}} ^{3}$ we have equations
\begin{equation}\label{eq_5_7}
(r_G^0(\lambda))^{*}(r_G^0(\lambda))^\top=\alpha^{2}\gamma'(\lambda)(F^{0}(\lambda )+G^{0}(\lambda))B_{1}^{\top}(F^{0}(\lambda)+G^{0}(\lambda)),
\end{equation}
\begin{equation} \label{eq_5_8}
\frac{1}{2 \pi} \int_{- \pi}^{ \pi} \left| \left \langle B_{1}, F^0( \lambda )-F_{1} ( \lambda ) \right \rangle \right|d\lambda= \delta,
\end{equation}
\begin{multline}\label{eq_5_9}
(r_F^0(\lambda))^{*}(r_F^0(\lambda))^\top=\\=
(\beta^{2} +\gamma'_{1} (\lambda )+\gamma'_{2} (\lambda
))(F^{0} (\lambda )+G^{0} (\lambda ))B_{2}^{\top}(F^{0} (\lambda)+G^{0} (\lambda )),
\end{multline}
where $\left| \gamma' ( \lambda ) \right| \le 1$ and
\[\gamma' ( \lambda )={ \mathrm{sign}}\; \left \langle B_{1} ,F^{0} ( \lambda )-F_{1} ( \lambda ) \right \rangle : \; \left \langle B_{1} ,F^{0} ( \lambda )-F_{1} ( \lambda ) \right \rangle \ne 0,\]
$\gamma'_{1}( \lambda )\le 0$ and $\gamma'_{1} ( \lambda )=0$ if $\langle B_{2},G^{0} ( \lambda \rangle > \langle B_{2},V( \lambda ) \rangle,$ $\gamma'_{2}( \lambda )\ge 0$ and $\gamma'_{2} ( \lambda )=0$ if $\langle
B_{2} ,G^{0} ( \lambda \rangle < \langle B_{2} ,U( \lambda ) \rangle.$

For the fourth pair $D_{1\delta}^{4}\times {D_{V}^{U}} ^{4}$ we have equations
\begin{equation}\label{eq_5_10}
(r_G^0(\lambda))^{*}(r_G^0(\lambda))^\top=(F^{0} (\lambda )+G^{0} (\lambda ))\left\{\alpha_{ij}\gamma_{ij}(\lambda ))\right\}_{i,j=1}^{T} (F^{0} (\lambda )+G^{0}(\lambda )),
\end{equation}
\begin{equation} \label{eq_5_11}
\frac{1}{2 \pi} \int_{- \pi}^{ \pi} \left|f^0_{ij}(\lambda)-f_{ij}^{1}( \lambda ) \right|d\lambda = \delta_{ij}, \; i,j=\overline{1,T},
\end{equation}
\begin{equation}\label{eq_5_12}
(r_F^0(\lambda))^{*}(r_F^0(\lambda))^\top=(F^{0} (\lambda )+G^{0} (\lambda ))(\vec{\beta}\cdot \vec{\beta}^{*}+\Gamma _{1} (\lambda )+\Gamma _{2} (\lambda ))(F^{0} (\lambda)+G^{0} (\lambda ))
\end{equation}
where $\left| \gamma_{ij} ( \lambda ) \right| \le 1$ and
\[\gamma_{ij} ( \lambda )= \frac{f_{ij}^{0} ( \lambda )-f_{ij}^{1} (\lambda )}{ \left|f_{ij}^{0} ( \lambda )-f_{ij}^{1}(\lambda) \right|} : \; f_{ij}^{0} ( \lambda )-f_{ij}^{1} ( \lambda ) \ne 0, \; i,j= \overline{1,T},\]
$\Gamma _{1} (\lambda )\le 0$ and $\Gamma _{1} (\lambda )=0$ if $G^{0}(\lambda )>V(\lambda ),$ $
\Gamma _{2} (\lambda )\ge 0$ and $\Gamma _{2} (\lambda )=0$ if $G^{0}(\lambda )<U(\lambda ).$

The following theorem and corollaries hold true.

\begin{theorem}{Theorem 4.}{}%
Let the minimality condition (\ref{minimal}) hold true. The least favorable spectral densities  $F^0(\lambda)$, $G^0(\lambda)$  in the classes $D_{1\delta} \times D_V^U$ for the optimal linear filtering of the functional $A\vec{\xi}$ are determined by relations
(\ref{eq_5_1}) -- (\ref{eq_5_3}) for the first pair $D_{1\delta}^{4}\times {D_{V}^{U}}^{1}$ of sets of admissible spectral densities;
(\ref{eq_5_4}) -- (\ref{eq_5_6}) for the second pair $D_{1\delta}^{4}\times {D_{V}^{U}}^{2}$ of sets of admissible spectral densities;
(\ref{eq_5_7}) -- (\ref{eq_5_9}) for the third pair $D_{1\delta}^{4}\times {D_{V}^{U}}^{3}$ of sets of admissible spectral densities;
(\ref{eq_5_10}) -- (\ref{eq_5_12}) for the fourth pair $D_{1\delta}^{4}\times {D_{V}^{U}}^{4}$of sets of admissible spectral densities;
constrained optimization problem (\ref{extrem}) and restrictions  on densities from the corresponding classes $D_{1\delta} \times D_V^U$.  The minimax-robust spectral characteristic of the optimal estimate of the functional $A\vec{\xi}$ is determined by the formula (\ref{4}).
\end{theorem}

\begin{corollary}{Corollary 6.}{}%
Assume that the spectral density $G(\lambda)$ is known. Let the function $F^0(\lambda)+G(\lambda)$ satisfy the minimality condition (\ref{minimal}). The spectral density $F^0(\lambda)$ is the least favorable in the classes $D_{1\delta}^k$, $k=\overline{1,4}$, for the optimal linear filtering of the functional $A\vec{\xi}$ if it satisfies relations (\ref{eq_5_1}) -- (\ref{eq_5_2}), (\ref{eq_5_4}) -- (\ref{eq_5_5}), (\ref{eq_5_7}) -- (\ref{eq_5_8}), (\ref{eq_5_10}) -- (\ref{eq_5_11}), respectively, and the pair $(F^0(\lambda), G(\lambda))$ is a solution of the optimization problem  (\ref{extrem}). The minimax-robust spectral characteristic of the optimal estimate of the functional $A\vec{\xi}$ is determined by formula (\ref{4}).
\end{corollary}

\begin{corollary}{Corollary 7.}{}%
Assume that the spectral density $F(\lambda)$ is known. Let the function $F(\lambda)+G^0(\lambda)$ satisfy the minimality condition (\ref{minimal}). The spectral density $G^0(\lambda)$ is the least favorable in the classes ${D_V^U}^k$, $k=\overline{1,4}$, for the optimal linear filtering of the functional $A\vec{\xi}$ if it satisfies relations (\ref{eq_5_3}), (\ref{eq_5_6}), (\ref{eq_5_9}), (\ref{eq_5_12}), respectively, and the pair $(F(\lambda), G^0(\lambda))$ is a solution of the optimization problem  (\ref{extrem}). The minimax-robust spectral characteristic of the optimal estimate of the functional $A\vec{\xi}$ is determined by formula (\ref{4}).
\end{corollary}

\section{Conclusions}
In the article we propose the methods of the mean-square optimal linear filtering of functionals which depend on the unknown values of a multidimensional stationary stochastic sequence based on observed data of the sequence with noise and missing values. In the case of spectral certainty where the spectral densities of the involved stationary sequences are known we derive formulas for calculating the spectral characteristic and the mean-square error of the optimal estimate of the functionals.
In the case of spectral uncertainty, where the spectral densities of the stationary sequences are not exactly known
while some special sets of admissible spectral densities are given, we apply the minimax-robust estimation method and derive relations which determine the least favourable spectral densities
and minimax-robust spectral characteristics.

\end{document}